,

# On the Simple Divisibility Restrictions by Polynomial Equation $a^n+b^n=c^n$ Itself in Fermat Last Theorem for Integer/Complex/Quaternion Triples

## Sandor Kristyan

*Research Centre for Natural Sciences, Eötvös Loránd Research Network (ELKH),*
*H-1117 Budapest, Magyar Tudósok körútja 2, Hungary*

Corresponding author: kristyan.sandor@ttk.mta.hu, kristyan.sandor@ttk.hu

**Abstract**. The divisibility restrictions in the famous equation $a^n+b^n=c^n$ in Fermat Last Theorem (FLT, 1637) is analyzed how it selects out many triples to be Fermat triple (i.e. solutions) if n>2, decreasing the cardinality of Fermat triples. In our analysis, the restriction on positive integer (PI) solutions ((a,b,c,n) up to the point when there is no more) is not along with restriction on power n∈PI as decreasing sets {PI }⊃{odd}⊃{primes}⊃{regular primes}, etc. as in the literature, but with respect to exclusion of more and more c∈PI as increasing sets {primes p}⊂{$p^k$}⊂{PI}. The divisibility and co-prime property in Fermat equation is analyzed in relation to exclusion of solutions, and the effect of simultaneous values of gcd(a,b,c), gcd(a+b,$c^n$), gcd(c-a,$b^n$) and gcd(c-b,$a^n$) on the decrease of cardinality of solutions is exhibited. Again, our derivation focuses mainly on the variable c rather than on variable n, oppositely to the literature in which the FLT is historically separated via the values of power n. Among the most famous are the known, about 2500 years old, existing Pythagorean triples (a,b,c,n=2) and the first milestones as the proved cases (of non-existence as n=3 by Gauss and later by Euler (1753) and n=4 by Fermat) less than 400 years ago. As it is known, Wiles has proved the FLT in 1995 in an abstract roundabout way. The n<0, n:=1/m, as well as complex and quaternion (a,b,c) cases focusing on Pythagoreans are commented. Odd powers FLT over quaternions breaks.

**Keywords.** Divisibility relations in Fermat Last Theorem; Cardinality of Fermat triples as a function of c; n<0 and n:=1/m Fermat cases; complex/quaternion Pythagoreans and Fermat triples; FLT breaks with odd powers over quaternions.

## INTRODUCTION

Fermat's Last Theorem (FLT) states [1-4] that no three positive integers (PI) a, b, and c satisfy the equation
$$a^n + b^n = c^n \qquad (1)$$
for any integer value of n greater than 2 in number theory. The cases n = 1 (elementary case, (1,1,2,n=1), (1,2,3,n=1), etc.) and n = 2 (Pythagorean triples, (3,4,5,n=2), (5,12,13,n=2), etc.) both have been known to have infinitely many solutions. (In both cases, n= 1 and 2, the number of solutions are countable-infinite, but using any threshold $c_{max}> 0$, there are more solutions for n= 1 than for n=2 if 1<c≤$c_{max}$.) Proofs by Fermat (n=4) and Gauss (n=3) and later for only small primes (>4) are about 1 page long and (elementary) algebraic. The general proof is necessary only for odd prime exponents (if $a^{2N}+b^{2N}=c^{2N}$ ⇒ $(a^2)^N+(b^2)^N=(c^2)^N$, and so on), as well as for co-primes (a,b,c) (via $(ka)^n+(kb)^n=(kc)^n$ wherein gcd≥k ⇒ $a^n+b^n=c^n$). By 1993, FLT had been proved for all primes less than 4·$10^6$ with the help of computers. The first successful but very long, lasting about 120 pages proof for any n>2 was formally published only in 1995 by Wiles [5] based on group theory. Over the centuries the opinion was born, that simply to juggle with algebraic powers (binomial theorem, product form of difference of powers, etc.) is useful in proof for cases with small n (like the known n=1-4 cases), but for all n>4 it is probably not enough, some other, more abstract ways have to be chosen. The famous statement from Fermat that, the short proof can be written to the margin of the page beside the equation seems to be a mistake, more, it is opinioned that looking for short poof of FLT is just looking for "gold of fools".

As indicated in the Abstract, our discussion is along the value c in Eq.1, not along with power n as in the literature:

**Definition 1.** When we talk about c in Eq.1, $c_{max}$ means that c=1,2,…, $c_{max}$ are discussed, and $c_0$ means that a particular c=$c_0$ value is discussed.

## SIMPLE SIEVES FOR FERMAT EQUATION WHICH EXCLUDES SOLUTIONS

Any zero (or ∞) element in set {a,b,c,n} ⊂ PI∪{0} leads to trivial solution for Eq.1. The a,b,c>0 and
$$0 < a < b < c < a+b \qquad (2)$$
holds among PI values. In Eqs.1-2, a=b can also occur for n=1, as well as "a" has been chosen smaller, but their symmetric appearance could fix the opposite b<a also, as well as the c<a+b restriction in Eq.2 for Eq.1 are evidenced in App.1. Textbook solutions for Eq.1 with 0<m≤M PI are (a,b,c,n)= (m,M,m+M,n=1) and (2mM,$M^2-m^2$,$M^2+m^2$,n=

2), in which for co-primes (a,b,c) further reductions may necessary. The trivial n=1 case has importance in the cardinality of solutions as a function of n, see Eq.8 below. All (m,M) PI pair links bijectively the (not necessarily co-prime) Fermat triples between n=1 and 2, so the cardinality is the same i.e. ∞, but at any fixed maximal c, there are less Pythagorean triples (n=2) than Fermat triples at n=1, see Table 1. The "smallest Fermat triple" (defined by the smallest "a" called $a_{min}(n)$ at constant n) (a,b,c,n=1 or 2) is (1,1,2,n=1) as 1+1=2, the next is (1,2,3,n=1), and Pythagorean triple (3,4,5,n=2) as the classical $3^2+4^2 = 5^2$ and the next is (5,12,13,n=2). Notice that $a_{min}(n=1)<a_{min}(n=2)$.

**TABLE 1**: Filtering the triples (a,b,c,n) at fixed $c=c_0=5$ in Eq.1 from being Fermat triples, i.e. it exhibits how the simple sieves work (Y= yes, it holds/filters). For $c=c_{max}=5$ all the (a,b)∈PI pairs are listed via Eq.2 as $0<a\leq b<c_{max}$. Last line (cumulative filter) leaves (a,b)= (2,4), (3,3) and (4,4) to be candidates for n>2, however, a<b if n>1, so only the (a,b)= (2,4) is not filtered. However, by Eq.6 the $n\leq c_0/2.164=2.31$, so the (a,b)=(2,4) is also filtered for n>2. Notice that this filtering does not need direct substitution for Eq.1, as well as the small $c_0$ value should not bother anyone, the filtering trend is same for any $c_0$ value.

| a | 1 | 1 | 1 | 1 | 2 | 2 | 2 | 3 | 3 | 4 |
|---|---|---|---|---|---|---|---|---|---|---|
| b | 1 | 2 | 3 | 4 | 2 | 3 | 4 | 3 | 4 | 4 |
| n=1 in Eq.1 holds as $a^1+b^1=c_0^1$ | | | | Y | | Y | | | | |
| n=2 in Eq.1 holds as $a^2+b^2=c_0^2$ (the Pythagoreans) | | | | | | | | | Y | |
| Eq.8 filters as n=1 and N=2 for n>2 | | | | Y | | Y | | | Y | |
| Eq.2 filters as $c_0\leq a+b$ | Y | Y | Y | | Y | | | | | |
| Cumulative filter by Eqs.2,8 | Y | Y | Y | Y | Y | Y | | | Y | |

Eqs.1-2 can be reformulated with x≡a/c and y≡b/c as

$$x^n + y^n = 1 \quad \text{and} \quad 0 < x < y < 1 \quad \text{wherein x and y are rational numbers.} \quad (3)$$

On the domain of real numbers, there are uncountable-infinite solutions (x,y) for Eq.3, however, now we restrict x and y to be rational numbers. If n= 1 or 2, Eq.3 has solutions that are purely (i.e. in both coordinates) rational (e.g. $(3/5)^2+(4/5)^2=1$) or purely irrational (e.g. at n=1 let x=1/π). The FLT (if true, but we know it is since 1995) yields that no pure rational (x,y) solution for Eq.3 if n>2. (Mixed solution, i.e. one of (x,y) is rational and the other is irrational does not exist by the definition of irrational numbers.)

**Lemma 1.** If (a,b,c,n>1) is Fermat triple, then c < a+b ≤ 2c-3, (if n=1 then c=a+b), more precisely

$$c < a+b \leq 2^{(n-1)/n}c = \{c\ (n=1),\ \sqrt{2}\ c\ (n=2), \ldots 2c\ (n=\infty)\} \leq 2c. \quad (4)$$

With x≡a/c and y≡b/c, the Fermat triples are in the strip with lines x+y=1 and $x+y=2^{(n-1)/n} \leq 2$.

The proof is as follows: Maximal candidates for "a" and "b" in Eqs.1-2 are a=c-2 and b=c-1 proving the cruder a+b= 2c-3. In Eq.3 x<y for n>1, as well as dividing the c < a+b ≤ 2c-3 by c yields 1 < x+y ≤ 2-3/c < 2. The finer value is detailed in ref.[6] and see App.1 for a corollary.

**Corollary of Lemma 1.** The reformulation of FLT with equivalent statements is as follows: 1.: The $c^nx^n+(b^n-c^n) = 0$ polynomial equation with PI (a,b,c,n) cannot have rational solution for n>2. 2.: Other algebraic form of Eq.1 is the $(x+k_1+k_2)^n-(x+k_1)^n-x^n = 0$ (or the $(x+k_1)^n-(x-k_2)^n-x^n = 0$, etc.) with $k_1,k_2 \in$PI having at least one real root if n>0 and all real roots are irrational if n>2, or the right hand side is never zero if x, $k_1, k_2 \in$PI and n>2. 3.: Another algebraic form of Eq.1 is by using c:=p+q, a:=p-q, p>q and p, q∈PI the $(p+q)^n = (p-q)^n+R^n$ holds, where $R^n\equiv \Sigma_{i=0}{}^n({}^n_i)p^i(q^{n-i}-(-q)^{n-i})$ by the binomial theorem, and ($R^n\in$PI but) R is always irrational if n>2. 4.: If 0<y<1 is rational and n>2 in Eq.3, then the $x = (1- y^n)^{1/n}$ is always irrational.

In relation to decide or filter if a number y is irrational in Eq.3, another Fermat theorem is useful which says that the only pairs of rational numbers that satisfy the equation $y^2=x^3-x=x(x^2-1)=x(x-1)(x+1)$ are (x,y)= (0,0), (1,0), (-1,0) – the proof is relatively simple (half page elementary algebra). What is important is that if x is real, then y is irrational, as well as it is used in the proof of FLT for n=4 in the literature. This curve is a special case of elliptic curves $y^2=x(x-A)(x+B)=x^3+A'x+B'$. The lengthy proof (completed by Wiles in 1995) is based on a special case of elliptic curves, the Frey's curve or Frey–Hellegouarch curve [7-9] (see App.2) associating hypothetical solutions (a,b,c,n) of Fermat's equation. This elliptic curve would have special properties due to the appearance of high powers of integers in its equation and the fact that Eq.1 would be an n-th power as well, importantly, any counter-example to FLT would probably also imply that an elliptic curve existed that was not modular. The modularity theorem states that elliptic curves over the field of rational numbers are related to modular forms. It was the crucial point and needs lengthy proof and the proof is far not simple. Wiles has proved the modularity theorem for semi-stable elliptic curves, which was enough to imply FLT. Any elliptic curve over rational numbers would have to be a modular elliptic curve, yet if a solution to Fermat's equation with non-zero (a,b,c,n>2) existed, the corresponding curve would not be modular,

,

resulting in a contradiction. Notice that, Fermat claimed to have a short (at least less than 1 page) proof, while modularity theorem as a part of proof is very complex, and incorporates the work of so many other specialists (e.g. refs.[5,7-10]) that it was suggested in prior to 1995. While the proof for low value n=3-4, etc. cases needed only classical algebra (and the high particular prime n values needed computation tricks), the general n with modular theorem needed the abstract group theory (see App.2).

To find how many rational points lie on a given algebraic curve comes from solving polynomial equations in two variables, however, one sees that the answer is not simple. There are easy examples of curves exhibiting the full range of behaviors, from having no rational points, to finitely many and infinitely many; it is far beyond the scope of current mathematics for an algorithm telling the exact number of rational points on any curve.

**Definition 2.** The cardinality of the set of solution triples (a,b,c,n) of Eqs.1-2 is denoted by non-negative arithmetic function $K(n) \geq 0$.

Since FLT has been proved back in 1995, we know that $K(n>2)=0$. FLT would have a shorter proof if monotonicity could be proved shortly as

**Conjecture 1.** $\qquad\qquad\qquad K(n) \geq K(n+1) \geq 0 \quad \text{for } n>0 \qquad\qquad\qquad$ (5)

As indicated above, the K(n) is countable-infinite for n=1 and 2 but no solution for n=3, i.e. K(3)=0 by Gauss or K(4)=0 by Fermat, and with these, the FLT follows immediately from Eq.5, as counterpart of the existing very long proof of FLT by modularity theorem [5,7-10]. It is also obvious, that with a $0<c_{max}<\infty$ threshold PI value, the number of solutions ($K(n,c_{max})$ is less if n=2 than if n=1, that is, $K(1,c_{max}) > K(2,c_{max})$, (see tables in ref.[6]). The hypothesis/conjecture is that the number of solutions is decreasing as n increases in Eq.1, that is, $K(n,c_{max}) \geq K(n+1,c_{max})$ for any given $c_{max}$, so $K(n) \geq K(n+1)$. Importantly, it has been well-known that, as the genus $((n-1)(n-2)/2)$ of curves gets higher, the number of rational points (or integral points) on arithmetic curves decreases [11]. It is written in many books about the arithmetic elliptic curves and especially in ref.[12], and exhibited in tables in ref.[6].

Eq.3 can be considered as a curve over real numbers, and we consider its restriction to rational numbers. An advantage of Eq.3 vs. Eq.1 in relation to filtering is that Eq.3 takes care on primitive triples, as well as it changes PI domain $(0,\infty)$ to rational domain $(0,1)$. Proving the decreasing property of K(n) (to prove FLT with it and not with e.g. the modularity theorem) Faltings's theorem (about 40 pages proof) has provided a milestone [13]. The Mordell conjecture (1922) says that a curve of genus greater than 1 over the field of rational numbers has only finitely many rational points, proved by G. Faltings (1983, 1984), and is now known as Faltings's theorem. Different proofs were found later, but yet the proofs are ineffective in that no bounds are given on the height of rational points or the upper bound for the number of rational points on the curve. In other words, Faltings's theorem tells only about finiteness of cardinality but not about the particular values of cardinality or how it changes with e.g. the genus. For Eq.3 it means that K(n) is finite if genus > 1. A corollary of Fallings's theorem [13] by Filaseta [14] is: For every n≥3 there exists a natural number m=m(n) (depending on n) such that if k≥m then FLT is true for kn. Using this, Granville [15] and Heath-Brown [16] established the theorem: For every real number x, let N(x) be the number of positive integers n≤x such that FLT is false for n, then $\lim_{x\to\infty} N(x)/x = 0$. These strongly support but do not proof Eq.5.

However, there are effective simple (i.e. the proof is simple and simple to use) sieves.

**Lemma 2.** If $a^n + b^n$ is an n-th power among the PI, then

$$n \leq a < b \quad \text{and} \quad n \leq c/\kappa \quad \text{with } \kappa \approx 2.164 \qquad (6)$$

The proof is simple and known in literature for κ=1. The finer κ value can be found in ref.[6] along with lemma:

**Lemma 3.** For any given $c_0 > 2$ integer in Eq.2, $(a,b,c_0,n>n_0)$ cannot be Fermat triple, where threshold $n_0$ is

$$1/\ln(2) \approx 1.443 \quad < \quad c_0/n_0 \approx 1.5/\ln(2) \approx 2.164 \quad < \quad 2/\ln(2) \approx 2.885 \qquad (7)$$

In case of n=3 (since Gauss) or generally n>2 (since FLT after 1995) Eq.6 states something on an empty set of PI, but Eq.6 is a useful filter when one works with K(n) to try prove FLT in another way. Eq.6 comes from the discrete nature of variables: If c is given, a large n enforces a=c-2 and b=c-1 to be maximal. For example, c:= 10, then at, let say, n=7 > 10/κ the $8^7+9^7= 6\,880\,121 < 10^7$ equality holds, and the irrational solution of $x^7+9^7=10^7$ is x≈9.11 > 8. We deal with PI, and Eq.6 supports (but does not prove) that at a threshold value $K(n,c_{max}) \geq K(n+1,c_{max})$, for which taking $c_{max} \to \infty$ the $K(n) \geq K(n+1)$. Obviously $0 \leq K(n,c_{max}) \leq K(n) \equiv \lim_{c_{max}\to\infty} K(n,c_{max})$. For example, $K(n=1,c_{max}=5)=2 > K(2,5)=1$, see Table 1.

**Lemma 4.** In the crude estimation $0 \leq K(n,c_{max}) << S(n,c_{max})$ the S is given in ref.[6] in which at least $S(n\to\infty, c_{max})=0$.

The following equality also supports the monotonic deceasing property of K(n), the proof is simple [6]:

**Lemma 5.** If Eq.2 holds ("≤" for n=1 an "<" for n>1) and 0 < n, N are PI, then

$$\text{If } n \neq N \text{ and } a^n + b^n = c^n \quad \Rightarrow \quad a^N + b^N \neq c^N \qquad (8)$$

Eq.8 is obvious among cases at n=1 and 2, Table 1 is a demonstration (full scan for $0 < a < b < c_0 \leq 5$ ). The n= 1 and 2 cases select out many candidates to be Fermat triple for n>2. (For example, $m^n+M^n \neq (m+M)^n$ if n>1 by Eq.8 or by binomial theorem, while equal if n=1, see App.3.) A conjecture is that

$$a_{min}(n) < a_{min}(n+1) \text{ if } n>0 \tag{9}$$

For example, $a_{min}(n=1)=1 < a_{min}(2)=3$, see Table 1. These very simple statements support, but does not prove Eq.5.

Eq.1 with n=2N $\in$ PI reads as $(a^N)^2+(b^N)^2=(c^N)^2$, and formula for Pythagorean triples provides that with all {m, M$\in$PI} the $\{a^N=M^2-m^2, b^N=2mM, c^N=M^2+m^2 \in PI\}$ are solutions, since $(a^N,b^N,c^N,2)$ is Pythagorean (Fermat) triple. If simultaneously happens that all three {a, b, c $\in$ PI or rational} also, then (a, b, c, 2N) is also a Fermat triple. It supports Eqs.5 and 9: If a triple happens to be Pythagorean, further filter is that $a=(M^2-m^2)^{1/N}$, etc. are also three integers, but since 1995 we know there are no such PI if N>1. If such a,b,c$\in$PI existed for an N>1, the alternate form $(a^2)^N+(b^2)^N= (c^2)^N$ would yield $(a^2,b^2,c^2,N>1)$ Fermat triple. Its corollary is that it is enough to consider FLT for odd n in Eq.1. Notice that, the 2N-power problem reduced to N-power conditions. The more general $(a^N)^{n_1}+(b^N)^{n_1}=(c^N)^{n_1}$ yields: If no PI solution for an $n_1$ (i.e. at least one in the brackets is irrational, e.g. $a^N$), there cannot be PI solution for $n=Nn_1$ ("a" must also be irrational), that is:

**Lemma 6**. It is enough to consider primes for n in Eq.1 for FLT (well known in literature).

## CERTAIN DIVISIBILITY PROPERTIES AS SIEVES FOR FERMAT EQUATION

As mentioned, it is enough to consider prime n powers (2, 3, 5, 7, 11,…, p, …) in Eq.1. Here, we provide additional relatively strong and simple selections or sieves to support FLT. It is based on the textbook relationships that, the

$$a^{2N+1}+b^{2N+1} = (a+b)(a^{2N}-a^{2N-1}b\pm\ldots-ab^{2N-1}+b^{2N}) \equiv (a+b)F(a,b,n=2N+1) \text{ for odd } n=2N+1 \in PI . \tag{10}$$

If Eq.1 holds for PI, the a+b and F(a,b,n) are simultaneous divisors as $(a+b)|c^n$ and $F|c^n$. The former provides a simple selection. Rearranging Eq.1 yields $b^n=c^n-a^n$, and the counterpart relationship is also useful: The

$$c^n-a^n = (c-a)(c^{n-1}+c^{n-2}a+\ldots+ca^{n-2}+a^{n-1}) \equiv (c-a)G(c,a,n) \text{ or } c^n-b^n= (c-b)G(c,b,n) \text{ for } n>1 \text{ and } n \in PI . \tag{11}$$

An expansion of Eq.1 comes from Eqs.10-11, interestingly for n=0, ±1, ±2, ±3, ±4,…, as

$$a^n+b^n = (a+b)(a^{n-1}+b^{n-1}) - ab(a^{n-2}+b^{n-2}) \tag{12}$$

$$c^n -a^n = (c -a)(c^{n-1}+a^{n-1}) + ca(c^{n-2} -a^{n-2}) \text{ and similarly for } a \to b . \tag{13}$$

The roots of polynomial equations F(a,x,n=2N+1)=0 are $az_i$ and of the G(x,a,n)=0 are -a and $az_i$ if n is even and only $az_i$ if n is odd with purely complex $z_i$ (i.e. irreducible over PI). (Notice that, unlike in binomial theorem, all coefficients are unity in F(a,b,n) and G(c,a,n).) Using Eq.6 and Lemma 1 the following holds:

**Theorem 1.** The constraint by divisor property is as follows: If a, b, c, n$\in$PI (let a<b) and Eq.1 hold, then c < a+b < 2c and

n ≤ a < b < c and n≤ c/κ, as well as

1.: If n=2N+1>2 (i.e. includes odd primes) then (a+b) and F(a,b,n) are simultaneous divisors as $(a+b)|c^n$ and $F|c^n$.

2.: If n>1 (i.e includes all primes) then two simultaneous divisor pairs hold as {$(c-a)|b^n$ and $G(c,a,n)|b^n$ simultaneously} and {$(c-b)|a^n$ and $G(c,b,n)|a^n$ simultaneously}.

3.: If n=2N then two simultaneous set of divisors hold as {$(c\pm a)$, $c^2-a^2$, G(c,a,n) and g(c,a,n) are divisors of $b^n$ simultaneously} and {same holds with interchanging a and b}.

4.: In points 2 and 3 the c-a>1 by Eq.2 (and cannot be the unit divisor), but the c-b=1 case provides a trivial unit divisor (e.g. in point 1 the $(a+c-1)|c^n$ and $F(a,c-1,2N+1)|c^n$ are enforced via Eq.1).

5.: Importantly, if one or more of the divisor properties in 1-3 are not satisfied, the (a,b,c,n) cannot be Fermat triple, i.e. FLT holds for it by this simple theorem.

Note: Theorem 1 is for n>2, and knowing that FLT is true, it states something on an empty set, like Eq.6, because there is no Fermat triple (a,b,c,n>2) solving Eq.1. However, it is interesting to check on a non-empty set, that is, when n=1 or 2. If n=1, i.e. a+b=c, then F=G=1 are trivial divisors and $(a+b)|c^1$, $(c-b)|a^1$ and $(c-a)|b^1$ trivially hold for Fermat triple (a,b,c=a+b,n=1). For Pythagorean triples n=2, i.e. $a^2+b^2=c^2$, then $a^2=c^2-b^2=(c-b)(c+b)$ and G=c+b, so $(c-b)|a^2$ and $G|a^2$, etc., e.g. $(c-a)|b^2$ and so on yielding that (a,b,$c^2=a^2+b^2$,n=2) is Fermat triple (seems trivial in this form), telling no detail how to obtain. That is, for example, (5,12,13,n=2) owns $13-12=1|5^2$, $13-5=8=2^3|12^2=(2^23)^2$, although (5+12)= 17DD$13^2$, but Theorem 1 for this latter asks odd n.

Elementary relationships between the parity (O$\equiv$ odd PI, E$\equiv$ even PI) of PI terms used in F, G and Eq.1 are:

**Lemma 7.** The parity restrictions in the Fermat equation in Eq.1, F(a,b,2N+1) and G(c,a or b,n) are as follows: Substituting even (E) or odd (O) PI values for a, b and c in Eq.1, F(a,b,2N+1) and G(c,a or b,n), the following relations hold (modulo 2, selecting the odd and even PI):

1.: If Eq.1 holds, then $O^n+E^n \equiv E^n+O^n \equiv O^n$ (i.e. cannot be $E^n$) or $O^n+O^n \equiv E^n$ (i.e. cannot be $O^n$) or $E^n+E^n \equiv E^n$, as well as in the latter case (a,b,c) are not co-primes, Eq.1 can be divided by $2^n$ at least (elementary).

2.: F(E,E,2N+1) $\equiv$ E, i.e. even a and even b yield even F values, F(E,O,2N+1) $\equiv$ F(O,E,2N+1) $\equiv$ F(O,O,2N+1) $\equiv$ O.

,

3.: $G(O,O,2N) \equiv E$,  $G(O,O,2N+1) \equiv O$,  $G(E,E,n) \equiv E$  and  $G(E,O,n) \equiv G(O,E,n) \equiv O$.

4.: $G(c,a,2N) = (c+a)g(c,a,2N)$  and  $G(c,a,4N+2)) = (c+a)F(c,a,2N+1)G(c,a,2N+1)$.

5.: In primitive (gcd(a,b,c)=1) Fermat triple (a,b,c,n=1) the c is even or odd (1+1=2, 1+2=3, etc.) and in primitive Pythagorean Fermat triple (a,b,c,n=2) the c is always odd simply from its generating equations (via the term 2Mm). However, what can we say about the parity of c when n>2, regardless FLT (i.e. such a "c" does not exist)?

As proof, the 1 and 2 are elementary. If a is odd and b is even, the $a^{2N}$, $a^{2N-k}b^k$ and $b^{2N}$ are odd, even and even resp., and F≡ odd + (2N-1)even + even ≡ odd; similarly if a is even and b is odd. If a and b are odd, the $a^{2N}$, $a^{2N-k}b^k$ and $b^{2N}$ are all odd, and F≡ odd + (2N-1)odd + odd ≡ (2N+1)odd≡ odd. The 3 is also elementary, via nO≡ E if n=2N and O if n=2N+1, nE≡E as well as (n-1)E+O≡ O+(n-1)E≡ O. Even more elementary, but crucial in the selections below that $O^n$=O, i.e. power of odd PI is odd, $E^n$=E for even PI, even number can have even and odd divisors (e.g. (6 and 9)|$2^4 3^5$), but odd number cannot have even divisor (e.g. 6DD$3^5$). The first part of 4 was reasoned above, the second comes from $(c-a)G(c,a,4N+2) = c^{4N+2}-a^{4N+2} = (c^{2N+1}+a^{2N+1})(c^{2N+1}-a^{2N+1}) = [(c+a)F(c,a,2N+1)][(c-a)G(c,a,2N+1)]$.

Notice that F|$c^n$ and G|($b^n$ or $a^n$) forecast the present of elliptic curves, since e.g. F|$c^n$ with y:=c and x:=a yield $y^n=(a+b)(x^{n-1}+...b^n)$. For example, if c=5 for Theorem 1, the Eqs.2 and 6 leave the residue candidate set for a+b as {2+4=6, 3+4=7}. The 3≤a selects and the residue is {a+b= 3+4=7}, however (3+4)DD$5^3$, so (a,b,c=5,n=3) is not a Fermat triple, a conclusion without direct calculation with Eq.1. By this theorem the divisor properties select further after Eqs.2 and 6 when K(n) is building up in Eq.5. If FLT is true (we know it is), then this theorem state divisor properties on an empty set of PI if n>2 (like Eq.6 with 2< n ≤ a, b), since there is no Fermat triples such those: If there were, then those would obey these divisor properties, and it is useful to use it for contradictions. Easy check for disobeying at least one of the divisor properties allows us to make a decision that a particular (a,b,c,n) cannot be Fermat triple. As a corollary, if n is odd and c is prime p, then (a+b)|$p^n$ should hold in Eq.1, however, divisor of $p^n$ in this case is {1, p, $p^2$, …, $p^n$}, but c=p < a+b < 2c=2p < $p^2$ also holds (Lemma 1), so a+b cannot be a divisor. More, if c is composite, especially c=$p^k$ with p prime and k∈PI, then (a+b)|$p^{kn}$ must hold, but the divisors of $p^{kn}$ are 1<p<$p^2$<…< c=$p^k$ < $p^{k+1}$ < … <$p^{nk}$. Again, c=$p^k$ < a+b < 2c=2$p^k$ < $p^{k+1}$ holds, that is, (a+b)|$c^n$ does not hold. Finally,

**Theorem 2.** The prime exclusions are as follows: For odd n=2N+1 and c=$p^k$ with p is prime and N, k∈PI simply, the (a+b)|$p^n$ cannot hold in Eq.1. An immediate consequence: (a,b,$p^k$,2N+1) cannot be a Fermat triple (as $a^{2N+1}+b^{2N+1}=(p^k)^{2N+1}$). Similarly and simply, by 0<c-a<b=$p^k$ and 0<c-b<a=$p^k$ and n>1 neither b (as $a^n+(p^k)^n=c^n$), neither a (as $(p^k)^n+b^n=c^n$) and neither both (as $(p^k)^n+((p')^k)^n=c^n$) can constitute Fermat triple. That is, FLT simply follows in these conditions.

As it is known, the (infinitely many) primes are very frequent in PI at the beginning, so it is a strong filter for K(n) in Eq.5. From the Prime Number Theorem, a modified [17-18] prime-counting function that gives the number of primes ≤ x for any real number x is $\pi(x) \approx x/\ln(x)+x/(\ln(x))^2$. For example in Table 1, at the $c_{max}$=5 the π(5)=3 i.e. the set {2,3,5}⊂PI and the 4=$2^2$ are annihilated by this theorem, i.e. all the candidates listed (see last line). Finally and shortly speaking, c≠$p^k$ in Eq.1 for n=2N+1>2 with N∈PI, p is prime and k∈PI, because from Eqs.1, 4 and 10 in $(a^n+b^n)/(a+b) = c^n/(a+b)$ the left is PI, but the right is non-PI because c<a+b<2c. It is a strong selection and support for FLT by a simple contradiction. For example $c_0$=3, $2^2$ and 5 are selected out in the small Table 1. We also mention that for every c there is a prime in [c,2c] (Chebyshev). Although half of the interval is even (i.e. non-prime), but e.g. in the small Table 1 the c=5, 7 are primes in [$c_{max}$=5, 2$c_{max}$=10].

**Theorem 3.** The co-prime exclusions, supporting FLT simply, are as follows: The followings hold for triples (a,b,c) solving Eq.1:

1.: If n=2N+1>2 and (a+b,c) is co-prime, then (a+b)|$c^n$ cannot hold in Eq.1, and as an immediate consequence, (a,b,c,n) cannot be a Fermat triple.

2.: If n>1 and (c-a,b) is co-prime, then (c-a)|$b^n$ cannot hold in Eq.1, and as an immediate consequence, (a,b,c,n) cannot be a Fermat triple. Now Eq.2 also provides that c-a >1.

3.: If n>1 and (c-b,a) is co-prime, then (c-b)|$a^n$ cannot hold in Eq.1, and as an immediate consequence, (a,b,c,n) cannot be a Fermat triple. Now c-b=1 is trivial case and must be investigated further.

4.: It is enough to consider co-primes (a,b,c) triples in Eq.1 (a known elementary fact, because if p is a common divisor, then Eq.1 can be divided by $p^n$). The 1-3 obviously include odd primes and simple consequences.

Theorem 3 comes from Theorem 1. Counter example: $60^2+91^2=109^2$, and indeed, (c-a,b)=($7^2$,7*13) and (c-b,a)= (18,60) are not co-primes so passes to Eq.1, and although (a+b,c)=(151,109) is co-prime, but n is even. If (a,b,c) does not have common divisor and c<a+b<2c, then it is very likely true for one of the doubles mentioned: For example, (3,4,5) is co-prime and (3+4,5), (5-3,4), (5-4,3) are yes, not, trivial co-primes, resp., etc.. Co-prime (a,b,c) means by definition that gcd(a,b,c)=1, i.e. not only pairwise. Important to note for the co-prime exclusion in Theorem 3 that, if gcd(a+b,c)=1 (analogously gcd(c-a,b)=1 or gcd(c-b,a)=1 except when c-b=1 trivially holds in the latter), then (a+b)DD$c^n$ (analogously (c-a)DD$b^n$ or (c-b)DD$a^n$) for sure, and it follows simply that FLT holds in these cases.

However, if gcd(a+b,c)≠1 (analogously gcd(c-a,b)≠1 or gcd(c-b,a)≠1), then still both, (a+b)|$c^n$ or (a+b)DD$c^n$ can hold, actual value makes the choice. (Similarly for the other two gcd, and see "Case c is even composite" below for a=7, b=20, c=21, n=5, (a+b)=27=$3^3$|$c^5$=$(3*7)^5$.) Theorem 1 uses simultaneous divisors, the (a+b)|$c^n$ and F|$c^n$ (analogously with G), and Theorem 3 uses gcd(a+b,c) values (and analogously other gcd values). To check (a+b)|$c^n$ and gcd(a+b,c) with other gcd mentioned is "much easier" than to check divisibility F| and G|. These allow to extend the ranges of triples (a,b,c) forbidden for Eq.1 and K(n).

The simplest composite case c=$p_1 p_2$ with primes $p_1$ and $p_2$ is as follows: If for example, c=3*7=21, then by Eq.6 n ≤ c/κ ≈ 9, let say n=5. The 36 divisors of $c^5$ are 1<3<7<9< c=3*7=21 < $3^3$=27 < $7^2$=49 <63...< $3^5 7^5$= 4084101. By Lemma 1 the 0 < a < b < c= 21 < a+b < 2c=42 and the constrain of Theorem 1 (point 1) is (a+b)|$c^5$. The residue set of these contain {a+b} = {22,23,...41}\{27}, because 27 is divisor. Eq.6 asks n=5≤ a in a+b along with a<b<c, so e.g. in this set 22= 5+17= 6+16= 7+15= 8+14= ...= 10+12, similarly for the others, and the last is 41= 5+36= ...= 20+21. These are not divisors of $c^5$ so FLT holds for them, i.e. the (a+b)|$c^{2N+1}$ is a strong selector. There is only one, the $p_1^3$=$3^3$=27 which is a divisor but not prime, it has to be consider further. It can be checked with direct individual calculations for Eq.1: For these FLT still holds, i.e. $a^5+b^5=21^5$ with a+b=27 has no PI solution, e.g. $5^5+22^5$= 5156757= $3^3$*31*61*101 ≠ 4084101= $(3*7)^5$= $21^5$, $7^5+20^5$= $3^3$*11*10831 ≠ $(3*7)^5$= $21^5$, etc.. Another way is via the simpler divisor check: The inspection of all possible cases of (a+b)=27= $3^3$|$c^n$= $(3*7)^5$ with all allowed (0<a<b<c=21) values, now $c^5$/(a+b)= $3^2 7^5$. Now gcd(a+b=27,c=21)= 3 > 1 i.e. not co-prime, and the check is in Table 2.

**TABLE 2**: Divisor check: Underlined number means reason for exclusion of co-prime (a,b,c=21,n) to be a Fermat triple because: 1.: non-PI, 2.: gcd(c-a,b)= gcd(c-b,a)= 1 i.e. division violation, 3.: gcd(a,b,c)≠ 1 i.e. not primitive triple (these values can form Fermat triples, but at a smaller c, which must be considered there and not in this Table).

| a | 7 | 8 | 9 | 10 | 11 | 12 | 13 |
|---|---|---|---|---|---|---|---|
| b=27-a | 20 | 19 | 18 | 17 | 16 | 15 | 14 |
| c-a=21-a | 14 | 13 | 12 | 11 | 10 | 9 | 8 |
| c-b=21-b | 1 | 2 | 3 | 4 | 5 | 6 | 7 |
| $a^5$/(21-b) | $7^5$ | $2^{14}$ | $3^9$ | $2^3 5^5$ | $\underline{11^5/5}$ | $2^9 3^4$ | $\underline{13^5/7}$ |
| $b^5$/(21-a) | $\underline{2^9 5^5/7}$ | $\underline{19^5/13}$ | $2^3 3^9$ | $\underline{17^5/11}$ | $\underline{2^{19}/5}$ | $3^3 5^5$ | $2^2 7^5$ |
| gcd(a,b,c=21) | 1 | 1 | $\underline{3}$ | 1 | 1 | $\underline{3}$ | 1 |
| gcd(c-a,b)= gcd(21-a,b) | 2 | $\underline{1}$ | 6 | $\underline{1}$ | 2 | 3 | 2 |
| gcd(c-b,a)= gcd(21-b,a) | $\underline{1}$ | 2 | 3 | 2 | $\underline{1}$ | 6 | $\underline{1}$ |

Now $O^5+E^5 ≡ E^5+O^5 ≡ O^5$ (mod 2) is satisfied for Eq.1, but (c-a)|$b^5$, (c-b)|$a^5$, F|$c^5$ and G(c, a or b,n)|($b^5$ or $a^5$ ) must also hold, as well as the a=5 and 6 drop as PI solution for Eq.1 by violating c>b. One can check one by one, e.g. the first c-a=21-5=16= $2^4$|$b^5$= $22^5$=$2^5 11^5$ but c=21<b=22 violating c>b, etc., and one needs to prove that every column has violation at a point (i.e. FLT holds). Column at c-b=1|$a^5$=$7^5$ is trivial, but c-a= 14= (2*7)DD$b^5$=$20^5$=$(2^2*5)^5$ is a violation, column at c-b= 2|$a^5$=$8^5$ holds, but c-a= 13DD$b^5$=$19^5$ is a violation, and so on. Finally column at c-b= 7DD$a^5$=$13^5$ is a violation although c-a=8= $2^3$|$b^5$=$14^5$=$(2*7)^5$ holds. Two columns are not selected out: the a=9 and 12. Here the (F)DD($c^5$) selects for example, as $c^5$/F(a,b,n=5)= (a+b)$c^5$/($a^5+b^5$)= 27$c^5$/($a^5$+$(27-a)^5$)= $7^5$/($3^3$11) if a=9 and $7^5$3/461 if a=12, i.e. non-PI values. For a general c this latter must be checked individually, but one can notice that a and b have alternating parity and alternating violations happen. The reason that (c-a) and (c-b) are strong filters or sieves as simultaneous divisors (of $b^5$ and $a^5$, resp.) is that parallel change one by one in a and b, so if divisor is ok in one, it is false in the other and vice versa. (All F(a,b,n=5) as well as all G(c=21,a,n=5) and G(c=21,b,n=5) are odd, so parity is not necessarily violated for division of odd $c^5$=$21^5$ and odd or even $b^5$ and $a^5$, but the actual values are not divisors, see numerical example for $7^5+20^5$ above.) Now violations were identified by individual calculations.

The total filter or cancellation comes from two co-prime properties, and it is not individual, but a general behavior: First, if (a,b,c) is not co-prime (gcd(a,b,c=21)>1), that case does not belong here, Eq.1 can be divided with (gcd(a,b,c=21))$^n$, and belongs to a smaller c'= c/gcd(a,b,c=21) < 21. Second, now (a+b=27,c=21) is not co-prime with gcd=3 (the reason for the table to inspect), but (c-a,b) and (c-b,a) are co-primes (altering if happens), and if one of it happens, for those (a,b,c) the Eq.1 does not hold (Theorem 3), that is, FLT holds. Finally, these two filters, {(a,b,c) is not co-prime (gcd>1)} and {(c-a,b) or (c-b,a) is co-prime (gcd=1)} together cancel all primitive (a,b,c=21) with a+b=27=$3^3$|$c^5$=$(3*7)^5$ in the table above to be Fermat triple for Eq.1. As a short inspection reveals, FLT holds for this residue set (a+b=27) too. In summary, checking the bottom three lines for the gcd values, every column has at least one underlined number (i.e. must be excluded). Notice that this co-prime test does not depend "strongly" on n, larger n increases the starting a=5 value yielding less triples for test, etc.. More, if n is lower, that also selects by e.g. for the 27=$3^3$|$(3*7)^n$ to be true, the lowest n is 3. One more consequence is that, if (a+b)|$c^n$ happens to be true (a frequent

situation as in the table above), the gcd(c-a,b) or gcd(c-b,a) excludes these (a,b,c,n>2 odd) values from Eq.1, so the $F(a,b,n=2N+1)DDc^n$ must hold which cannot be proved simply (such as $y^{2N+1}=kF(a,x,2N+1)$ has no rational points). The above simple and general behavior of gcd yields the FLT to be true without direct calculation of Eq.1.

One case in the above table must be commented, the column (a,b,c,n)= (7, 20, 21,5). Now, $(a+b)=27=3^3|c^5=(3*7)^5$ and $(c-b)=1|a^5=7^5$ hold, but the latter by trivial reason. However, $(c-a)=14=(2*7)DDb^5=(2^25)^5$ and it yields that FLT holds for this (a,b,c); although gcd(c-a,b)=gcd(2*7,$2^2$*5)=2≠1 and gcd(a,b,c)=(7,$2^2$5,3*7)=1 do not select. Generally, this is the case $a^n+(c-1)^n=c^n$ and a>1 by App.1, so a=2,3,4,…,c-2<b=c-1. The b=c-1 and c has opposite parity, so Lemma 7 by parity restrictions is not violated if "a" is odd. Now a=7, or generally a=3,5,7,…c-2 now c-2=19, but selector c<a+b<2c selects further as a=3,5,7,9,11,13, finally, the too small "a" excludes. A little more detail is as follows. The a=c-2 is too large since $(c-2)^n+(c-1)^n=?=c^n$ is same as considering the difference $(C-1)^n+C^n–(C+1)^n= C^n+\Sigma_{k=0}^{n-1}(^n_k)C^k>0$, i.e. non-zero by adding positive terms, so (c-2,c-1,c,n>2 odd) does not solve Eq.1. Furthermore, without generality, the $x^n+(c-1)^n=c^n$ with n=3,5,7,…/κ, now n=3,5,7,9 and c=21 has the (irrational) x≈ 10.804, 15.4632, 17.590, 18.420, resp., so a=7 is too small if b=20. This equation has the approximate solution if c is large (from about 1000) $x^n= c^n-(c^n-nc^{n-1}±…(-1)^n)≈ nc^{n-1}$, i.e. x≈ $n^{1/n}c^{(n-1)/n}$. For n=3, c/2, ∞ this is $3^{1/3}c^{2/3}≈ 1,4422c^{2/3}$, $c/2^{2/c}≈c$, c, resp., and in the table above (and generally) a=7 << [x= $5^{1/5}21^{4/5}$]$_{integer\ part}$=15, i.e. cannot be a Fermat triple.

Cases with even composite c are as follows: If e.g. c=2*3*7=42 is even composite, then by Eq.6 n ≤ c/κ≈ 19, let say n=7. There are many divisors of $c^7$, e.g. 1,2,3…12,… < c=42 <…< 49=$7^2$ <…< 64=$2^6$ <… <2c=84 <…$42^7$, let us pick the odd a+b=49 and even a+b=64 from interval (c,2c). The constrain (a+b)|$c^7$ holds for these and must be considered further (with F|$c^7$ or with the easier gcd). For (a+b)∈(c,2c) where (a+b)DD$c^7$, the FLT simply holds. Eq.6 asks n=7≤ a in a+b along with a<b<c=42, so e.g. in these two cases, the 49= 8+41= 9+40= 10+39= … = 22+27= 23+26= 24+25, and 64= 23+41= 24+40= 25+39= … = 29+35 =30+34= 31+33. The 49 and 64 are divisors and not primes, it has to be considered further. It can be checked with direct individual calculations for Eq.1 that, for these FLT still holds. However, simpler and rather general (i.e. not individual) check with divisors (i.e. it can be generalized to any c) allows to exclude all these in agreement with FLT. The inspection of the two series with (a+b)= (49 or 64)|$c^n$= $(2*3*7)^7$, wherein $c^7/(a+b)= (2*3*7)^7/7^2=(2*3)^77^5$ and $(2*3*7)^7/2^6= 2(3*7)^7$ is as follows. The gcd(a+b=49,c=42)= 7 > 1 and gcd(a+b=64,c=42)= 2 > 1, i.e. not co-primes, so Theorem 3, right now, cannot exclude immediately from Eq.1. In the sum of a+b=49, the parity of a and b oppositely altering, and violates the $O^7+E^7 ≡ E^7+O^7 ≡ O^7$ (mod 2) rule, since c=42 is even. In the a+b=64 case the parity of a and b is the same, they agree with $O^7+O^7 ≡ E^7+E^7 ≡ E^7$ (mod 2), however, when a and b are even, Eq.1 can be divided by at least $2^7$ now, since c=42 is even. Finally, it is enough to consider odd a and b in a+b=64. Notice that these exclusions are initiated by the (even) c but not its particular value, being very systematic and general as the next step.

For case c=64 the check is in Table 3.

**TABLE 3**: Divisor check: Underlined number means reason for exclusion of co-prime (a,b,c=42,n) to be a Fermat triple: 1.: gcd(c-a,b)= gcd(c-b,a)= 1 i.e. division violation, 2.: gcd(a,b,c)≠ 1 i.e. not primitive triple (non-primitive values can be Fermat triple, but at a smaller c which must be considered there and not in this Table).

| a | 23 | 25 | 27 | 29 | 31 |
|---|---|---|---|---|---|
| b=64-a | 41 | 39 | 37 | 35 | 33 |
| c-a=42-a | 19 | 17 | 15 | 13 | 11 |
| c-b=42-b | 1 | 3 | 5 | 7 | 9 |
| gcd(a,b,c=42) | 1 | 1 | 1 | 1 | 1 |
| gcd(c-a,b)= gcd(42-a,b) | <u>1</u> | <u>1</u> | <u>1</u> | <u>1</u> | 11 |
| gcd(c-b,a)= gcd(42-b,a) | <u>1</u> | <u>1</u> | <u>1</u> | <u>1</u> | <u>1</u> |

All are excluded by gcd(c-a,b)=1 or gcd(c-b,a)=1, although the many prime numbers in this region also selects in gcd. Notice that in fact, these do not depend on n, the n can have all its allowed values, n=1,3,5,7,…,19. Altogether, FLT holds for (a,b,c,n) in this range by this simple consideration and no need to substitute (a,b,c,n) into Eq.1. Again, these two filters, {(a,b,c) is not co-prime (gcd>1)} and {(c-a,b) or (c-b,a) are co-primes (gcd=1)} together cancel all primitive (a,b,c=42) with a+b= (49 or 64)|$42^7$ to be Fermat triple for Eq.1 by the simple and general consideration above. As this short inspection reveals, FLT holds. Notice that in this simple test (generalization is trivial) we have proved that $a^n+b^n=42^n$ does not have PI solution if n=2N+1=3,5,7,9,…,19 up to the large number $42^{19}≈ 6.85*10^{30}$ by the last one. The first column gcd(c-b,a)= gcd(1,23)=1 is the case $a^n+(c-1)^n=c^n$ as above in "simplest composite case $c=p_1p_2$" (detailed there), where the gcd(1,k)=1 trivially holds if k∈PI.

,

Important is that, in these two cases considered we started with a c>3 PI value, and in the allowed (trivial) range {0<a<b<c < a+b < 2c and n≤a,b,c/2.164} 1., If $c=p^k$ with prime p, then Theorem 2 simply proves the FLT, 2., If c is composite, then the gcd(a,b,c)=k>1 reduces the triple (a,b,c)~(a/k,b/k,c/k), that is, that has to be considered at a smaller c value, 3., If k=1 and if $(a+b)DDc^n$, then FLT simply follows by Theorem 3, as well as if $(a+b)|c^n$, then one or both, gcd(c-a,b) and gcd(c-b,a) is unity, and the Theorem 3 also excludes this (a,b,c,n) triple from Eq.1 to be a Fermat triple, so FLT follows simply again. The latter two gcd constitute a filter "comb" with pairwise {≠1,=1,…,≠1,=1} together with {=1,≠1,…,=1,≠1} or both as {=1,=1,…=1,=1}. (The latter is similar to the trivial fact wherein one starts anywhere in PI, the odd/even values are obtained alternatively in any direction.) Only the b=c-1 constitutes distinct but general case, however, in FLT the a≠1 and/or b≠c-1 if n>2 by simple algebra. Finally, Theorem 1 makes the constrain on Eq.1 with n>2 primes that $(a+b)|c^n$, $F|c^n$, $(c-b)|a^n$, $G(c,b,n)|a^n$, $(c-a)|b^n$ and $G(c,a,n)|b^n$ must simultaneously hold which cannot happen, at least one is violated. More, violation happens already among the algebraically "easier to handle" $(a+b)|c^n$, $(c-b)|a^n$ and $(c-a)|b^n$ as a result of the simultaneous behavior of gcd(a,b,c)=1, gcd(c-a,b) and gcd(c-b,a).

We summarize the divisor restriction properties. The divisor properties, more exactly constrains $((a+b)$ and $F)|c^n$, etc. have been exhibited via examples (with composites c=3*7 and 2*3*7), however it is clear that, a general proof (for odd and even c) is just the same supporting FLT via a much shorter and easier way wherein the parities of a, b and c have to be discussed only for divisibility. Shortly: Eq.1 can always be decomposed with polynomials as $c^n = a^n+b^n = (a+b)F(a,b,n=2N+1>2)$, $b^n = c^n-a^n = (c-a)G(a,b,n>1)$ and similarly for $a^n$ forcing $gcd(a+b,c^n)>1$, $gcd(c-a,b^n)>1$ and $gcd(c-b,a^n)>1$ simultaneously. If (a,b,c) is co-prime triple for Eq.1 (which is enough to consider), then at least one of (a+b,c), (c-a,b) and (c-b,a) is co-prime double (along with 0 < a < b < c < a+b < 2c and n ≤ a, b, c/2.164), so one of the first terms in the doubles (or pairs) cannot divide $c^n$, $b^n$ and $a^n$, resp., violating Theorem 3 (with gcd>1) for E.1, that is, FLT holds. (The n=2N+1>2 in F includes all odd primes enough for proving FLT.) This divisor property or constrain is responsible that Eq.1 has no PI solution for n>2. Eq.1 has no PI solution if n=3, (Gauss or reasoning above), but the equations $a^3+b^3=kc^3$, $A^3+B^3+C^3=D^3$, $A^3+B^3=C^3+D^3$ or $A^3+B^3=C^2$ have PI solutions, since the extra variable (k and D) counterbalances the restrictions or constrains by divisibility: For example, with low values

$$3^3+5^3 = 19*2^3 \quad \text{and} \quad 3^3+4^3+5^3 = 6^3 \quad \text{and} \quad 1^3+12^3 = 9^3+10^3 = 1729 \quad \text{and} \quad 4^3+8^3 = 24^2. \qquad (14)$$

In Eq.14, indeed $19*2^3/(3+5)= 19$, $(6^3-5^3)/(3+4)= 13$, $(3^3+4^3)/(6-5)= 91$, $1729/(1+12)= 133$, $1729/(9+10)= 91$, $4+8= 12|24^2$, i.e. divisors by Eq.10, but Eq.2 is meaningless, (or in 4th the 4+8<24 i.e. Eq.2 formally fails but $4+8 > (4^3+8^3)^{1/3} ≈ 8.32$). The increasing genus makes monotonic decrease in cardinality of solutions as K(n)≥K(n+1)≥0 in Eq.1, while the divisibility constrains cause step-wise declining in cardinality as K(1)=K(2)=∞ and K(n>2)=0 (see tables, App.4).

Eqs.5 and 9 were called conjectures above, but knowing that FLT is true (Wiles 1995), those are theorems. The internal relationships between these equations are as follows. If we can prove Eq.5 or 9 in "shorter" ways, these yield FLT. Eq.5 is the stronger, however Eq.9 is also strong: For n= 1 and 2 the Eqs.5, 6 and 9 are true as it can be verified via e.g. Pythagoreans. Eq.6 simply provides $n<a_{min}(n)$ but does not prove Eqs.5 and 9. Eq.9 does not prove Eq.5 and vice versa, but a full proof of Eq.9 may be as simple. However, both, Eq.5 or 9 yields FLT because from n=3 the K(3)=0 (Gauss), that is $a_{min}(3)=∞$. Eq.5 is in agreement with that, the number of rational solutions (cardinality) decreases with increasing genus. Theorems 1-3 makes simple selections to support FLT.

## EXTENSIONS TO n<0, n:=1/m, AND COMPLEX OR QUATERNION (a,b,c) DOMAINS

For a given a, b, n ∈PI, Eq.1 can always be solved for real c, but for <u>integer (or rational) c</u> the FLT governs. The condition n=0, ±1, ±2, ±3… for Eqs.12-13 indicates that FLT can be <u>extended for range n→-n</u> in powers as well. Dividing Eq.1 by $(abc)^n$ yields $(ac)^{-n}+(bc)^{-n}=(ab)^{-n}$, i.e. if (a,b,c,n) is Fermat triple, then (ac,bc,ab,-n) also, furthermore, if gcd(a,b,c)=1, then gcd(ac,bc,ab)=1 also. On the other side, if $a^{-n}+b^{-n}=c^{-n}$, then multiplying by $(abc)^n$ yields $c(a^n+b^n)^{1/n}=ab$, i.e. there must be a c" such as (a,b,c",n) is Fermat triple. Similarly, $x^n+y^n=(xy)^n$ or $x^{-n}+y^{-n}=1$, and one has to consider if it has pure rational solution, etc.. For example, $(1/2)^4+(1/5)^4=(1/c)^4$ yields irrational c ≈ 1.9874. If FLT holds (we know it does), it holds for n→-n also. For example, (3,4,7,n=1) is Fermat triple, and (ca,cb,ab,-n)= (21,28,12,n=-1) also, or (3,4,5,n=2) is Fermat triple, and (15,20,12,n=-2) is also. Along with App.1:

**Lemma 8.** The extension of FLT and the cardinality of solutions for n<0 are as follows: If (a,b,c,n) with gcd(a,b,c)=1 is relative prime (primitive) Fermat triple, then (a'=ca,b'=cb,c'=ab,-n) with gcd(a',b',c')=1 is also a primitive Fermat triple, (for n=-1,-2 see more explicit form in ref.[2], here simpler). The cardinality at a given n is the same for –n up to primitive property. If |n|>2, then no solution for Eq.1 by FLT, that is, not only for n>2 but for n<2 too. The role of a' and b' are symmetric, but instead of Eq.2, 0<c'<a'≤b' holds for n=1 and 0<c'<a'<b' for n=2, furthermore, 0<a<a'=ac, etc.

When we allow the exponent n:=1/m for some integer m, we have the "inverse Fermat equation" $A^{1/m}+B^{1/m}=C^{1/m}$. In the case in which the $m^{th}$ roots are required to be real and positive, all solutions are reported in ref.[19] as $A=rs^m$, $B=rt^m$ and $C=r(s+t)^m$ with r,s,t∈PI and gcd(s,t)=1. Indeed, $A^{1/m}+B^{1/m}= sr^{1/m}+tr^{1/m}= (s+t)r^{1/m}= C^{1/m}$. For example, if

r=s=t=1 and m=2, then indeed $\sqrt{1}+\sqrt{1}=\sqrt{4}=2$, if r=s=m=2 and t=1, then indeed $\sqrt{8}+\sqrt{2}=\sqrt{18}$. Notice that square roots are involved what contains irrational number behind. If $r=R^{km}$, then the equation in fact is $sR^k+tR^k=(s+t)R^k$. Seemingly trivial, but important is the "all solution" property if we make back "derivation". To compare with Eq.1, $sR^k+tR^k=a^n+b^n$ (for the sum and not one by one) and $(s+t)R^k=c^n$. The latter enforces c=R=s+t and k=n-1. Substituting this into the left of Eq.1 yield $sc^{n-1} + (c-s)c^{n-1} = c^n$ with s=1,2,3,,,c-1, in fact this latter holds by itself. However, one cannot separate in this way for n-powers over PI as $a^n=sc^{n-1}$ and $b^n=(c-s)^{n-1}$, e.g. in $a=s^{1/n}c^{(n-1)/n}$ the s can be an n-power, but $c^{(n-1)/n}$ will be irrational except when n=1, then $c^{(n-1)/n}=1$, as well as (n-1)DDn except at n=2 (Pythagorean power).

Eq.1 on integer-complex numbers (A≡ $a+a_1i$ over PI with i≡$\sqrt{(-1)}$, etc.) and integer-quaternions (A≡ $a+a_1i+a_2j+a_3k$ over PI, etc.) domains is commented in ref.[6]; trivial is when $a_2=a_3=0$ reduces the quaternions to its complex sub-space. Literature of FLT in relation to complex numbers is not new but not long [2]. According to the derivation for n=2, the complex-Pythagoreans can be generated in the same way, but not (in the same way) for the quaternions, because the latter is non-commutative, see details in App.5. The $(a+a_1i)i = -a_1+ai$ makes sense (i.e. not-trivial) that if $(z_1,z_2,z_3,n)$ is a complex Fermat triple, than $(i^kz_1,i^kz_2,i^kz_3,n)$ also with integer k, even though FLT restricts for n>2.

**Lemma 9.** The extension of FLT to complex numbers is as follows: Solutions for Eq.1 over integer-complex numbers have the same form as over PI if n=1 or 2.

App.5 contains the proof, but examples are useful. Extending the domains in Lemma 9 if n=1 to negative a and b and complex A and B values are trivial, e.g. 2-1=1, falling into the PI case as 1+1=2, however the n=2 case has some curiosities. M=1 and m=i yield the simplest complex-Pythagorean as (a,ai,0,n=2) with integer a. The m=1 and M=1+i integer-complex numbers with positive components generate triple (-1+2i, 2+2i, 1+2i, n=2), i.e. negative component occurs. The m=1 and M=2+i generate reducible (2+4i, 4+2i, 4+4i, n=2) to (1+2i, 2+i, 2+2i, n=2):

$$(1+2i)^2 + (2+i)^2 = (2+2i)^2. \qquad (15)$$

Eq.15 has positive components, as well as notice that PI values m=1 and M=2 generate the "smallest" (3,4,5,n=2) so integer-complex number can generate triples with smaller coefficients in components than this one, recall Eq.9 for the same PI domain, now with different (PI vs. complex) domains. It comes from more components (a vs. a and $a_1$, etc.). Eq.1 with integer-complex domain as $A^2+B^2=C^2$ breaks into real and imaginary components (with the help of binomial theorem) and

**Lemma 10.** A Diophantine equation system of two variables is generated by a complex FLT equation. Eq.1 with integer-complex numbers generate the following Diophantine equation system of two

$$\Sigma_{k=0}^n \binom{n}{k}(a^{n-k}(a_1i)^k + b^{n-k}(b_1i)^k - c^{n-k}(c_1i)^k) = 0 \qquad (16)$$

$$a^2 + b^2 + c_1^2 = a_1^2 + b_1^2 + c^2 \quad \text{and} \quad aa_1 + bb_1 = cc_1 \qquad \text{by n=2} \quad (17)$$

$$a^3 + b^3 + 3cc_1^2 = 3aa_1^2 + 3bb_1^2 + c^3 \quad \text{and} \quad 3a^2a_1 + 3b^2b_1 + c_1^3 = a_1^3 + b_1^3 + 3c^2c_1 \qquad \text{by n=3} \quad (18)$$

where a, $a_1$, b, $b_1$, c, $c_1$ ∈ PI and i is the complex imaginary unit. The case n=2 in Eq.17 can be solved via integer-complex-Pythagoreans in App.5, for example the above calculated m=1 and M=2+i generate solution 1,2,2,1,2,2.

Notice that, the left of Eq.17 as $a^2+b^2-c^2 = a_1^2+b_1^2-c_1^2$ is satisfied by any two Pythagorean over PI, but the right side not necessarily holds, e.g. (3,4,5,n=2) with (5,12,13,n=2) the left is zero, but the right side is $aa_1+bb_1- cc_1 = 3*5+4*12-5*13= -2\ne 0$ or $3*12+4*5-5*13= -9\ne 0$. With the example in Lemma 10, $a^2+b^2-c^2= a_1^2+b_1^2-c_1^2 = 1\ne 0$ i.e. not-Pythagoreans over PI, but satisfies both in Eq.17, i.e. together Pythagoreans over integer-complex numbers. FLT holds for Eq.1, but if it holds for integer-complex numbers, it must also be proved somehow or extend. The left side in Eqs.17-18 includes Eq.1, but the more terms allow more flexibility, that is, existence of solutions if n>2 may be possible, or at least the threshold separating the solutions/no-solutions may be higher than n=2. (Compare this to Eq.14 versus FLT.) We do not discuss it further, but we mention that for Eq.18, for example, the low values $(1+3i)^3 + (-2+2i)^3 - (1+2i)^3 = 1$ (zero imaginary part), i.e. 1,3,-2,2,1,2 solves right of Eq.18. The $(2+i)^3+(-2+2i)^3-(3+i)^3 = i$ (zero real part), i.e. 2,1,-2,2,3,1 solves left of Eq.18. Both sides of Eq.18 contains extra terms with respect to Eq.1, allowing low value solutions at least for one of it. These can be found simply e.g. by scanning low values. However, solution must be found for Eq.18 which satisfy both equations simultaneously. The $\{a_1= b_1= c_1= 0\}$ or $\{a= b= c= 0\}$ reduce Eqs.16-18 to Eq.1 over PI wherein FLT holds, a triviality, but the latter still says a little fact: FLT is true for pure imaginary integer-complex numbers (it simply comes also from multiplying Eq.1 by $i^n= \pm1$ or $\pm i$). We do not discuss the important FLT for integer-complex numbers case further for generalization.

In relation to Pythagoreans (n=2) we can state the following properties: Substituting generating formula over PI for (M,m) and ($M_1,m_1$) yields PI-Pythagorean (a,b,c,n=2) in which a and b can be interchanged and PI-Pythagorean ($a_1,b_1,c_1,n=2$) in which $a_1$ and $b_1$ can also be interchanged. The left of Eq.17 is satisfied by 0=0, the right yields $aa_1+bb_1-cc_1 = -2(Mm_1-mM_1)^2= 0$, that is, $Mm_1=mM_1$ must hold for $(a+a_1i, b+b_1i, c+c_1i,n=2)$ to be complex-Pythagorean. But this condition means that this complex-Pythagorean is (a+kai, b+kbi, c+kci, n=2) with k∈PI, the first can be primitive, e.g. (3+3ki, 4+4ki, 5+5ki, n=2). Another way is, if $a_1=ka$, $b_1=kb$ and $c_1=kc$ in the right of Eq.17, then (a,b,c,n=2) is Pythagorean and the left also holds. For the counterpart, $ab_1+ba_1-cc_1 = -M^2(M_1-m_1)^2 -m^2(M_1+m_1)^2 + 2Mm(M_1-$ $m_1)(M_1+m_1)= 0$ enforcing $M_1 \pm m_1=0$ at the same time which yields $M_1=m_1=0$ trivial case, $(a+a_1i, 0, a+a_1i, n=2)$. Finally,

**Lemma 11.** Simple properties between real and imaginary parts in complex Pythagoreans follow: With integer-complex numbers $a+a_1i$, $b+b_1i$ and $c+c_1i$ (i.e. $a, a_1, b, b_1, c, c_1 \in$ PI), for the complex-Pythagorean (n=2) solutions of Eq.1 the following properties hold (neglecting trivial cases).

1.: If (a,b,c,n=2) PI-Pythagorean, then (qa,qb,qc,n=2) is PI/ complex/ quaternion-Pythagorean with q is PI/ complex/ quaternion (trivial; the non-commutative property of quaternions does not come up in powers).

2.: If (a,b,c,n=2) and $(a_1,b_1,c_1,n=2)$ are both PI-Pythagoreans, then dot/scalar product of two associated vectors $(a,b,c)(a_1,b_1,-c_1) = aa_1+bb_1-cc_1 = 0$ if $a_1=ka$, $b_1=kb$ and $c_1=kc$ (essentially the same PI-Pythagoreans) and $\neq 0$ otherwise (essentially different PI-Pythagoreans; see tables of Pythagoreans on the web).

3.: If $(a+a_1i, b+b_1i, c+c_1i, n=2)$ is complex-Pythagorean, then, 1., (a,b,c,n=2) and $(a_1,b_1,c_1,n=2)$ are both PI-Pythagoreans (enforced by left of Eq.17) but k∈PI exists such as $a_1=ka$, $b_1=kb$ and $c_1=kc$ (enforced by right of Eq.17), i.e. essentially the same two PI-Pythagoreans involved as (a+kai, b+kbi, c+kci, n=2), as well as one of it (the first) can be primitive, 2., Nor (a,b,c,n=2) and neither $(a_1,b_1,c_1,n=2)$ are PI-Pythagoreans (by left of Eq.17), 3., no case when one is Pythagorean and other is not between (a,b,c,n=2) and $(a_1,b_1,c_1,n=2)$ (by left of Eq.17).

**Lemma 12.** The symmetric properties of complex solutions are as follows: The symmetric property in relation to Eq.1 vs. Eqs.16-18 extends from A↔B i.e. "if $(a, a_1, b, b_1, c, c_1)$ solves Eq.16-18, then $(b, b_1, a, a_1, c, c_1)$ also" to that, "if $(a, a_1, b, b_1, c, c_1)$ solves Eq.16-18, then $(a_1, a, b_1, b, c_1, c)$ also". The latter is that, if complex triple $(z_1, z_2, z_3, n)$ solves Eq.1 then $(z_1', z_2', z_3', n)$ also, where $z_i'$ is obtained by interchanging its real and imaginary parts (mirror to line z=u+ui).

For example, with m=1 and M=1+i the $(2+2i)^2 + (-1+2i)^2 = (1+2i)^2$ from Pythagorean formula (App.5), then parallel by this lemma the $(2+2i)^2 + (2-i)^2 = (2+i)^2$ also holds. This lemma is true hypothetically also, i.e. when solution does not exist at larger n values. The purely imaginary integer complex triples $(a_1i,b_1i,c_1i,n)$ solving Eq.1 reduce to integer triples $(a_1,b_1,c_1,n)$ via dividing by $i^n$, so FLT holds, but see the opposite behavior of quaternions in Theorem 4 below.

App.5 explains that generating equations for PI-Pythagoreans can be used simply for complex-Pythagoreans, and with non-commutative restrictions for quaternion-Pythagoreans. A simple case is exhibited next: The $2m^2M^2+2M^2m^2= (2mM)^2$ for the subset of quaternions in which commutative property holds. For example, subset containing q= a+bi+bj+bk≡ (a,b,b,b) elements with a,b∈PI the $q_1q_2=q_2q_1$ holds and preserves the form (a,b,b,b), i.e. stays in the same subset (commutative subset). For example, with $m=(m_1,1,1,1)$ and $M=(M_1,1,1,1)$ the mM= Mm= $(m_1M_1-3, m_1+M_1, m_1+M_1, m_1+M_1,)$; the components are in commutative PI domain generating commutative mM quaternion products. If $M_1=2$ and $m_1=1$, then $M^2 \pm m^2= (2,1,1,1)^2 \pm (1,1,1,1)^2= (1,4,4,4) \pm (-2,2,2,2)= (-1,6,6,6)$ or (3,2,2,2), 2mM= 2(2,1,1,1)(1,1,1,1)= (-2,6,6,6) and indeed

$$(3+2i+2j+2k)^2 + (-2+6i+6j+6k)^2 = (-1+6i+6j+6k)^2 \tag{19}$$

because $(3+2i+2j+2k)^2 = (-3,12,12,12)$, $(-2+6i+6j+6k)^2 = (-104,-24,-24,-24)$ and $(-1+6i+6j+6k)^2 = (-107,-12,-12,-12)$. Another case is when M=2m along with m and M are quaternions yielding $2m^2M^2+2M^2m^2= (2m)^4$. With these $M^2 \pm m^2= 5m^2$ or $3m^2$ and 2mM= $4m^2$, but this case falls into $(3m^2)^2+(4m^2)^2=(5m^2)^2$ i.e. point 1 of Lemma 11. Interesting and it can be proved in elementary way:

**Lemma 13.** The Pythagorean imaginary part in quaternions is as follows: If integer-quaternion a+bi+cj+dk contains integer-Pythagorean (b,c,d,n=2), e.g. (-3,4,5,n=2) in its imaginary/vector part, then any products or $n^{th}$ power of these preserve the PI-Pythagorean property in the imaginary part, but destroy the primitive property (if any) when a'+a"≠1 as

$$(a'+bi+cj+dk)(a''+bi+cj+dk) = (a'a''-2d^2) + (a'+a'')bi + (a'+a'')cj + (a'+a'')dk . \tag{20}$$

Eq.20 yields the statements in this Lemma for the products, $n^{th}$ powers and primitiveness by recursive application.

These kind in Eq.20 yield $(a+bi+cj+dk)(-a+bi+cj+dk) = -(a^2+2d^2) \in$ PI purely real and (a+bi+cj+dk)+(-a+bi+cj+dk) = 2bi+2cj+2dk purely imaginary quaternion. For example, $(1+3i+4j+5k)^8 = 2923601 +(3s)i+(4s)j+(5s)k$ with s= -862792; because now |a|=1 < |d|=5, the |s| ∈$[d^8= 390625, (2d^2)^4= 6250000]$. Right of Eq.17 contains the same kind of terms appearing in quaternion product yielding Lemma 13. Lemma 13 is useful to consider FLT in certain subsets of quaternions, e.g. for Eq.1 the $(a +a_1i +a_2j +a_3k)^n = t +(sa_1)i +(sa_2)j +(sa_3)k$ eliminates the power if the vector part is PI-Pythagorean. The (d,0,d,n=2) and (0,d,d,n=2) are trivial PI-Pythagoreans, but yield non-trivial case in Eq.20. If (b,c,d,n=2) is PI-Pythagorean, than the (c,b,d,n=2) also, allowing permutation between b and c in Eq.20. Even more interesting that b, c and d can be permuted in the imaginary parts of Eq.20 but $a'a''-2d^2$ is invariant for the permutation (d the largest number in Pythagorean used); as well as the same permutation must be used everywhere in Eq.20. For purely imaginary values, Eq.20 yields $(0+bi+cj+dk)^3=-2d^2(0+bi+cj+dk)$ i.e. power only stretches the vector, and $(0+bi+cj+dk)^4=4d^4(1+0i+0j+0k)$.

Using both, vector and coordinate notation for quaternions, for purely imaginary quaternions the $(0+bi+cj+dk)^{2N}=$
$=((.)^2)^N = (-(b^2+c^2+d^2),0,0,0)^N = (-(b^2+c^2+d^2))^N \in PI$ with $N= 1,2,\ldots$ holds. Eq.1 for these yields $(0,b_1,c_1,d_1)^{2N} + (0,b_2,c_2,d_2)^{2N} = (0,b_3,c_2,d_3)^{2N}$ between quaternions, and it reduces to $(b_1^2+c_1^2+d_1^2)^N + (b_2^2+c_2^2+d_2^2)^N = (b_3^2+c_3^2+d_3^2)^N$ between PI. The latter is subcase of a, b, and c in Eq.1 on PI where FLT holds. If $N=1$, it is not a strong restriction, i.e. when $n= 2N= 2$. There are many quaternion-Pythagoreans like this. For example, if $N=1$, then $(0,1,0,0)^2+(0,0,1,0)^2 = (-2,0,0,0) = (0,0,1,1)^2$ holds; compare to other quaternion forms in Eq.19 or up to Eq.15 in sub-space. More,

$$(0+bi+0j+0k)^2 + (0+0i+bj+0k)^2 = (0+0i+bj+bk)^2 \text{ and } b \in \{integers\}. \quad (21)$$

Extension of this lemma from even n to odd n (unlike for FLT over PI) is not as simple, because $(0+bi+cj+dk)^{2N+1}= (-(b^2+c^2+d^2))^N(0+bi+cj+dk)$. Substituting this into Eq.1 yields $b_1(b_1^2+c_1^2+d_1^2)^N + b_2(b_2^2+c_2^2+d_2^2)^N = b_3(b_3^2+c_3^2+d_3^2)^N$ via comparing the real and imaginary parts for i and analogously for j and k. However, it is not the form in Eq.1 over PI, further elaboration is necessary to extend the FLT for integer-quaternions and odd n. For proving existence of solution, simple computer search enough for low value solutions, for example, if $N=1$, then $(0,-1,-1,0)$, $(0,0,1,-1)$ and $(0,-1,0,-1)$ solves. The cyclic multiplication rule in quaternions allows to break FLT for quaternions. Finally,

**Theorem 4.** The FLT for purely imaginary quaternions follows: For Eq.1 with "purely imaginary integer-quaternions" and even $n= 2N \in PI$, the FLT holds as no quaternion (of this kind) solution if $N>1$. Importantly, there are "purely imaginary integer-quaternions" that solves Eq.1 for odd $n= 2N+1 \in PI$, i.e. FLT does not hold generally for "set of quaternions":

$$(0-bi-bj+0k)^{2N+1} + (0+0i+bj-bk)^{2N+1} = (0-bi+0j-bk)^{2N+1} \text{ and } b \in \{integers\}. \quad (22)$$

More, simultaneous permutation in imaginary parts also yields equality in Eq.22.

We do not discuss the important general case of quaternion-Pythagoreans when the non-commutative property is fully present, nor the FLT for integer-quaternions generally. Notice that, in complex ($A=a+a_1i$, B, C, $n=2$) and quaternion ($q_1$, $q_2$, $q_3$, $n=2$) cases of Pythagoreans, the generating PI coordinates (m, $M \in PI$) yield integer components ($|a|$, $|a_1| \in PI$, etc.), i.e. not only PI components. More, Eq.9 is in effect also in comparison between (3,4,5,n=2), Eq.15 and Eqs.21-22.

## CONCLUSIONS

Simple restrictions, coming out instantly from FLT equation itself by divisors via Eqs.10-11, in the search for Fermat triples from Diophantine $a^n+b^n=c^n$ simply yield that FLT holds if any or two or all of a, b and c are powers of primes ($p^k$), as well as if (a+b,c) or (c-a,b) or (c-b,a) are co-primes. Comments have been made on PI/ complex/ quaternion-Pythagoreans as well as on the relationship between FLT and quaternions.

## APPENDIX

**Appendix 1.** a.: There are many trivial and elementary cases included/associated with FLT we do not deal with. For example, if $a=0$ then all $b=c \in PI$ solves Eq.1, if $n=0$ then $a^0+b^0=2 \neq c^0=1$ in Eq.1, if n=even and (a,b,c,n) is Fermat triple, then ($\pm a$, $\pm b$, $\pm c$, n) is also a Fermat triple (this is a property on empty set if $n>2$ by FLT), etc..

b.: If $a=b$ in Eq.1, than $2a^n=c^n$ and $c=2^{1/n}a$, the latter yields irrational c for $n>1$ and integer for $n=1$. That is, $a<b$ (or alternatively $b<a$) if $n>1$, but $a \leq b$ (or alternatively $b \leq a$) if $n=1$ in Eq.1.

c.: The $a+b>c$ or $\leq c$. If $a+b>c$ then $(a+b)^n= a^n+b^n +(na^{n-1}b+\ldots) >c^n$ by binomial theorem, however decreasing (minoring) the left one cannot simply conclude that $a^n+b^n=c^n$ happens or not, that must be investigated further, and that is a historically difficult task. If $a+b \leq c$, then $(a+b)^n= a^n+b^n +(na^{n-1}b+\ldots) \leq c^n$ with binomial theorem, and $a^n+b^n < c^n$ by simply minoring the left knowing that $0<a<b$, which means now the important $a^n+b^n \neq c^n$. Alternatively, if $a+b \leq c$, then $c^{n-1}a+c^{n-1}b \leq c^n$, minoring the left with $a<b<c$ yields $a^n+b^n<c^n$, i.e. $a^n+b^n \neq c^n$, that is, FLT holds in this case.

If $a,b \in PI$ is given, Lemma 1 or Eq.4 provides $c < a+b \leq 2c-3$, and Eq.2 forces c to be between minimal $b+1$ and maximal $a+b-1$. Finally, $c \in [b+1, a+b-1]$ along with Eq.6, see App.3 also. Fermat triples exist if $n=2$ called Pythagorean triples, and indeed there are triples with $c= b+1$, e.g. (3,4,5,n=2) or (19,180,181,n=2). Furthermore, if $c \leq c_{max}=100$, there are $K(2,c_{max}=100)= 16$ triples, among them there are 6 with $c=b+1$, that is 37.5%. If $c_{max}=300$, $K(2,c_{max}=300)= 47$, among them there are 9 with $c=b+1$, that is 19.15 %, a decreasing density.

**Appendix 2.** In view of the divisibility properties and theorems for FLT equation in this work (using elementary algebraic tools), Frey's curve or Frey–Hellegouarch curve [7-9], $y^2=x(x-a^n)(x+b^n)$, associating hypothetical solutions (a,b,c,n) of FLT equation does not exists.

The Taniyama-Shimura conjecture, since its proof (1995) known as the modularity theorem, says formally that, for every elliptic curve $y^2=Ax^3+Bx^2+Cx+D$ over the rational numbers, there exist non-constant modular functions $f(z)$ and $g(z)$ of the "same level" such that $[f(z)]^2=A[g(z)]^2+Cg(z)+D$. The most characteristic property of modular function in its definition is that, it is a holomorphic complex-valued function, $f(z)$, on the upper half complex plane with $f((az+b)/(cz+d))= (cz+d)^k f(z)$ and $k \in PI$; for odd k the $f=0$, as well as $f(z+1)=f(z)$ is an example.

,

**Appendix 3.** If a=1 with a given c, Eq.2 forces b up to b=c-1, however it is a Fermat triple by Eq.1 as $1^1+(c-1)^1 = c^1$ so by Eq.8 the $1^n+(c-1)^n = 1+(c-1)^n \neq c^n$ for n>1, finally a>1 if n>1 (which is elementary). By Eq.2, the allowed c up to this point is c=4=$2^2$, 5, 6, 7…, and by Theorem 2 these reduced to c= 6, 7,… .

**Appendix 4.** If c is prime and 0<a≤b<c, then (a,b,c) is co-prime. For example, if c=11 prime and n=5 odd, divisor of $11^5$ are {1,11,121,1331,$11^4$,$11^5$}, but c=11 < a+b < 2c=22, so (a+b) does not divide $11^5$ meaning that (a,b,c=11,n=5) is not a Fermat triple. Notice that for divisors of $p^n$, only the $1<p<p^2$ has to be considered further, i.e. the "low divisors". In relation to larger numbers, if e.g. c=1013 prime is considered, then simply 1013 < a+b < 2026 by Eq.2 and 2 < n ≤ a < b < c=1013 as well as primes n < 1013/κ≈ 1013/2.164≈ 468 by Eq.6 can be candidates for (a,b,1013,n) in Eq.1 only to check FLT; direct calculations to check other numbers if Eq.1 holds are waste of time. Notice that these cases contain numbers large as $c^n= (1013)^{468}$≈ $4.2*10^{1406}$ in magnitude, even FORTAN programs for direct handling of these numbers in direct check for Eq.1 is not simple. For a composite larger number and odd n=2N+1, for example, the $(p-b)^{2N+1} + b^{2N+1} = (11^{1934}13^{1955})^{2N+1}$ has no PI solution for p>b and N if 11≠p≠13 prime, because prime a+b= (p)DD$(11^{1934}13^{1955})^{2N+1}$, as well as Eq.6 also selects as no PI solution if $2N+1 > 11^{1934}13^{1955}/κ$. However if e.g. p=$11^{1926}13^{1959}$ there could be if N>0, but those will be excluded by F(p-b,b,2N+1) or by divisors c-a, c-b, etc.. Finally, with restriction c= $p^k$ with prime p in Eq.1, FLT holds in simple way (Theorems 1-3) which can be extended to composite values for c, see examples above.

**Appendix 5.** If n=1, the solutions for Eq.1 are (a,b,c,n=1)= (m,M,m+M,n=1) among PI, integer-complex numbers and integer-quaternions (trivial). The Pythagoreans generated along with $(M^2+m^2)^2-(M^2-m^2)^2= 2m^2M^2+2M^2m^2$. Notice that, squares ($M^2$ and $m^2$) are necessary for a, b and c in Eq.1 if n=2 versus $1^{st}$ powers when n=1; supporting Eq.9 between n=1 and 2. The right is $4m^2M^2=(2mM)^2$ for PI and integer-complex numbers yielding the form in Eq.1 and any m and M in PI or integer-complex numbers can be used. However, integer-quaternions are (commutative in $m^2$ and $M^2$, but) non-commutative generally in multiplication, i.e. $m^2M^2 \neq M^2m^2$, so one has to solve/elaborate the $2m^2M^2+2M^2m^2= q^{2\ or\ 4}$ further to get the form in Eq.1 which restricts the allowed integer-quaternions for m and M.

Furthermore, quadratic 2-2, 3-2 and 4-1 (not necessary primitive) equations can be obtained via manipulating the sums of $\pm(\pm r \pm t \pm s)^2$ followed by r:=$r^2$, etc.. The r=s=t=q=1 yields $2^2=1^2+1^2+1^2+1^2$, more, irrational r=√2, s=1/√2, p=2√2 and q=3/√2 yield $5^2+2^2+4^2+6^2=9^2$, generally for PI (and r=√L with L∈PI, etc.) and integer-complex variables:

$$(r^2+s^2-t^2)^2 + (2rt)^2 = (r^2-s^2+t^2)^2 + (2rs)^2 \text{ and } (-r^2+2s^2+t)^2 + (r^2-2s^2+t)^2 + (4rs)^2 = (r^2+2s^2+t)^2 + (r^2+2s^2-t)^2 \quad (23)$$

$$(r^2-s^2-t^2-q^2)^2 + (2rs)^2 + (2rt)^2 + (2rq)^2 = (r^2+s^2+t^2+q^2)^2. \quad (24)$$

## ACKNOWLEDGMENTS

Financial and emotional support for this research from OTKA-K-NINCS 2015-115733 and 2016-119358 is kindly acknowledged. Thanks to Szeger Hermin for her help in typing the manuscript. The subject has been published in the American Institute of Physics Conference Proceedings, ICNAAM 2022.

## ABBREVIATIONS

DD= does not divide (e.g. a|b= "a divides b" vs. aDDb= "a does not divide b"); FLT= Fermat's Last Theorem;
Fermat triple= (a,b,c,n) containing PI and satisfying Eq.1 (although a 4 variable vector); gcd= greatest common divisor(s);
K(n)= number of solutions (a,b,c) at fixed n if any for Eq.1; O= odd PI and E= even PI;
PI= positive integer(s), we avoid notation N or $N_0$ for this because e.g. Eq.8 or 10;
PI/integer/complex/quaternion-Pythagorean= all variables in Eq.1 with n=2 are PI/integer/integer-complex/integer-quaternion.

## REFERENCES

1. S. Singh: Fermat's Last Theorem, Fourth Estate (Publisher), 1997
2. https://en.wikipedia.org and https://mathworld.wolfram.com
3. L. E. Dickson: History of the Theory of Numbers, Volume 2., New York, Chelsea, 1952
4. G. H. Hardy, E. M. Wright: An Introduction to the Theory of Numbers, Oxford, 1965
5. A. Wiles: Modular elliptic curves and Fermat's Last Theorem, Annals of Mathematics, **141** (3), 443–551 (1995)
6. S. Kristyan: Note on the inductive proof for FLT, Part I-II.: ICNAAM-2021, Am. Inst. Phys. Conf. Proc. (2023)
7. Y. Hellegouarch: Polska Akademia Nauk. Instytut Matematyczny. Acta Arithmetica, **26** (3), 253–263 (1974)
8. G. Frey: J. Reine Angew. Math., **331**, 185–191 (1982)
9. G. Frey: Annales Universitatis Saraviensis. Series Mathematicae, **1** (1), (1986)
10. K. Ribet: Annales de la faculté des sciences de Toulouse Sér. **5**, 11 (1), 116–139 (1990)
11. Y. Manin: Rational points on algebraic curves, Russian Mathematical Surveys, **19** (6), 75-78 (1964)
12. J. H. Silverman: Advanced Topics in the Arithmetic of Elliptic Curves, Springer-Verlag, 2009
13. G. Faltings: Inventiones Mathematicae, **73** (3), 349-366 (1983) with erratum **75** (2), 381 (1984)
14. M. Filaseta: An application of Faltings's results to FLT, C.R. Math. Rep. Acad. Sci. Canada, **6**, 31-32 (1984)

,

15. A. Granville: C.R. Math. Rep. Acad. Sci. Canada, **7**, 55-60 (1985)
16. D. R. Heath-Brown: FLT for "almost all" exponents, Bull. London Math. Soc., **17**, 15-16 (1985)
17. S. Kristyan: Prime Numbers Are Not Erratic, Am. Inst. Phys. Conf. Proc., **1863**, 560013 (2017)
18. S. Kristyan: Note on the cardinality of primes and twin primes, Am. Inst. Phys. Conf. Proc. **1978**, 470064 (2018)
19. M. Newman: A radical Diophantine equation, Journal of Number Theory, **13** (4), 495–498 (1981)